\documentclass[12pt, a4paper]{article}
\usepackage{enumerate, amsmath, amssymb, amscd, paralist,array,url,amsthm}
\usepackage[all,web,poly,curve,frame]{xy}
\newcommand{\R}{\mathbb{R}}
\newcommand{\Z}{\mathbb{Z}}
\newcommand{\N}{\mathbb{N}}
\newcommand{\Q}{\mathbb{Q}}
\newcommand{\C}{\mathbb{C}}

\newcommand{\B}{\mathcal{B}}

\newcommand{\Ll}{\mathcal{L}}

\newcommand{\Hom}{\operatorname{Hom}}
\newcommand{\Pic}{\operatorname{Pic}}
\newcommand{\Oo}{\mathcal{O}}
\begin{document}
\title{Cohomology of line bundles on a toric variety\\ and constructible sheaves on its polytope}
\author{Nathan Broomhead}
\date{}
\maketitle
\begin{abstract} We explain a method for calculating the cohomology of line bundles on a toric variety in terms of the cohomology of certain constructible sheaves on the polytope. We show its effective use by means of some examples.
\end{abstract}
\section{Introduction}
Let $X$ be a toric variety (not necessarily smooth) over an algebraically closed field of characteristic zero, with $T$ the embedded torus and let $D$ be a torus-invariant Cartier divisor on $X$. It is a well known result (see \cite{Dani, Fulton}) that the cohomology $H^p(X,\Oo(D))$ splits into a direct sum of weight spaces indexed by the character lattice $M = \Hom(T,\C^*)$,
\begin{equation}\label{dsd} H^p(X,\Oo(D)) \cong \bigoplus_{m \in M}H^p(X,\Oo(D))_m
\end{equation} 
There is a theorem of Demazure (see \cite{Dani}) which says that each weight space can be written as a local cohomology group, calculated on the fan of $X$,
$$\quad H^p(X,\Oo(D))_m \cong H^p_{Z(\psi,m)}(| \Delta |, \C)$$
where $Z(\psi,m):= \{ v \in | \Delta | \mid \langle m , v \rangle \geq \psi(v) \}$ and $\psi = \psi_D$ as defined in section~\ref{Defs}.

In this article we show that when $X$ is projective, each weight space can also be written in terms of the cohomology of certain constructible sheaves on the polytope $P_X$ corresponding to $X$. 
\newtheorem{main}{Theorem}[section]
\begin{main}\label{maint} For all $p \geq 0$ there are canonical isomorphisms:
\begin{equation} \label{mthm}
H^p(X,\Oo(D))_m \cong H^p(P_X,j_! \underline{\C}_{W})
\end{equation}
where $W=W(m,D)$ is the complement of $Z(m,D)$, a union  of closed maximal dimensional faces of $P_X$, $j$ is the open inclusion of $W(m,D)$ into $P_X$ and $ \underline{\C}_{W}$ is the constant sheaf on $W$.
\end{main}

Furthermore there is a short exact sequence of constructible sheaves (see \cite{Hartshorne,Iversen}),
\begin{equation} \label{seseq3} 0 \longrightarrow {j}_!{j}^*\underline{\C}_{P_X} \longrightarrow \underline{\C}_{P_X} \longrightarrow {i}_*{i}^*\underline{\C}_{P_X} \longrightarrow 0
\end{equation} 
where $i$ is the inclusion of $Z(m,D)$ into $P_X$, which induces a long exact sequence of cohomology. Since ${j}^*\underline{\C}_{P_X} = \underline{\C}_{W} $ and ${i}^*\underline{\C}_{P_X} = \underline{\C}_{Z} $, we can see that the long exact sequence is:
\begin{equation*} \label{long} \dots \longrightarrow H^{p-1}(Z(m,D), \C) \longrightarrow H^{p}(P_X,j_! \underline{\C}_{W}) \longrightarrow H^{p}(P_X, \C) \longrightarrow \dots
\end{equation*}
whence we can calculate $H^p(P_X,j_! \underline{\C}_{W})$ in terms of the complex cohomology on $Z(m,D)$. 

Finally we remark that although the above has been stated in the projective case, if $Y$ is a quasi-projective toric variety with a given embedding as an open subset of some projective toric variety $X$, then there is a corresponding open subset $P_Y$ of the polytope $P_X$ and, restricting to $Y$ and $P_Y$ respectively, the proof works as in the projective case.

{ \bf Acknowledgement}: I would like to express my deepest thanks to my supervisors: Dr. Gregory K. Sankaran who constructed the example that is considered in section~\ref{example}; and Dr. Alastair King, for many long discussions and much valuable input.
 
\section{Definitions and notation}\label{Defs}
In this section we define the objects and notation used in the proof.
Let $X = X(\Delta)$ be an n-dimensional projective toric variety corresponding to some complete fan $\Delta$ in a lattice $N \cong \Z^n$. We write $\Delta(r)$ for the set of $r$-dimensional cones in $\Delta$, and label the set of generators in $N$ of the $1$-dimensional cones by $\{e_i \mid i \in I \}$. It is well known (see \cite{Fulton}) that there is a 1-1 correspondence between prime torus invariant divisors and the elements of $\Delta(1)$, whence we shall denote these divisors $\{E_i \mid i \in I\}$ respectively. Let $ M := \operatorname{Hom}(N, \mathbb Z)  \cong \mathbb Z^n $ be the dual lattice to $N$, with dual pairing $ \langle \cdot , \cdot \rangle $, and let $M_\mathbb R := M \otimes_\mathbb Z \mathbb R \cong \mathbb R^n$.
 
Choose a divisor $A = \sum_{i \in I} a_i E_i$ which is Cartier and ample.
There is a polytope in $M_\mathbb R$ associated to $A$, given by:
\begin{align*}
P_A :&= \{u \in M_\mathbb R \mid \langle u , e_i \rangle \geq -a_i \; \forall i\}\\
    &= \{u \in M_\mathbb R \mid  u \geq \psi_{A} \; \text{on} \; | \Delta |= N_{\R}\}
\end{align*}
where the $\Delta$-linear support function $\psi_A : | \Delta | \rightarrow \R$  is determined by the property
$ \psi_A(e_i) = -a_i  $ for all $i \in I$.
For any cone $\sigma \in \Delta$ we define:
\begin{align*}
{T_\sigma} :&= \{u \in M_\mathbb R \mid \langle u , e_i \rangle = -a_i \; \forall e_i \in |\sigma|, \; \text{and} \; \langle u , e_i \rangle > -a_i \; \text{otherwise} \}\\
    &= \{u \in M_\mathbb R \mid  u = \psi_{A} \; \text{on} \; | \sigma |, \; u > \psi_{A} \; \text{on} \; |\Delta|\backslash | \sigma | \}
\end{align*}
This is open in its closure $\overline{T_\sigma} = \{u \in M_\mathbb R \mid  u = \psi_{A} \; \text{on} \; | \sigma |, \; u \geq \psi_{A} \; \text{on} \; |\Delta|\backslash | \sigma | \} $ and is thus locally closed. It can be seen that this defines a stratification $\{T_\sigma\}_{\sigma \in \Delta}$ of $P_A$ by locally closed sets indexed by the cones $\sigma \in \Delta$. Since $A$ is ample, each of these $T_\sigma$ is non-empty.
There is a natural partial order on the strata given by:
$$ T_\tau \leq T_\sigma \iff T_\sigma \subseteq \overline{T_\tau} $$
\newtheorem{poslem}{Lemma}[section]\label{Poslemma}
\begin{poslem}
$$ T_\tau \leq T_\sigma \iff \tau \subseteq \sigma $$
\end{poslem}
\begin{proof}
By definition $\tau \subseteq \sigma$ trivially implies ${T_\sigma} \subseteq \overline{T_\tau}$.
Conversely suppose $u \in {T_\sigma} \subseteq \overline{T_\tau}$. Then 
$  u > \psi_{A} $ on $ |\Delta|\backslash | \sigma | $ and also $ u = \psi_{A} $ on $ | \tau |$
whence $| \tau | \subseteq | \sigma |$. \end{proof}
A very similar argument shows that $\overline{T_\sigma}$ can be written as a union of the strata:
\begin{equation} \label{clloc} \overline{T_\sigma} = \bigcup_{\tau \supseteq \sigma}T_\tau \end{equation}
Note: As a topological object with a stratification by locally closed subsets, $P_A$ is actually independent of the choice of ample divisor $A$ on $X$ and depends only on $\Delta(X)$ (ie. on $X$). Thus in general we write $P_X$ for this object when we don't want to specify a particular ample divisor and we treat the strata $\{T_\sigma\}_{\sigma \in \Delta}$ as well defined subsets of this. From equation~(\ref{clloc}) we see that any closed maximal dimensional face $\overline{T_{\langle e \rangle}}$ is a union of strata depending only on the fan and is therefore also well defined in $P_X$. We denote it by $F_e \subset P_X$.
\newtheorem{cover}[poslem]{Lemma}
\begin{cover}\label{ocover}
There exists a cover of $P_X$ by open sets of the form $ V_\sigma := \bigcup_{\tau \subseteq \sigma}T_\tau$
\end{cover}
The proof of this is straightforward and left to the reader.
\newtheorem{contractible}[poslem]{Lemma}
\begin{contractible}\label{cntct}
For any $\sigma \in \Delta$, $V_\sigma$ is a contractible space.
\end{contractible}
\begin{proof} $V_{\sigma}$ is convex: Let $a,b \in V_{\sigma}$ so $a \in T_{\tau}$ and $b \in T_{\gamma}$ for some $\tau,\gamma \subseteq \sigma$. Let $e \in |\Delta|$, then
$$ \langle ta+(1-t)b , e \rangle = t\langle a , e \rangle + (1-t)\langle b , e \rangle \geq \psi_A(e) $$
for $t \in (0,1)$, since $\langle a , e \rangle \geq \psi_A(e)$, and $\langle b , e \rangle \geq \psi_A(e)$, with equality if and only if $e \in | \tau \cap \gamma|$. Thus $ta+(1-t)b \in T_{\tau \cap \gamma} \subseteq V_{\sigma}$. \end{proof}

Let $D = \sum_{i \in I} d_i E_i$ be a torus invariant divisor which is Cartier and fix some $m \in M$. We define a (closed) subset of $P_X$ by:
$$Z(m,D):= \bigcup_{ \{i \in I \mid \langle m , e_i \rangle < -d_i \} } F_{e_i}$$
a union  of closed maximal dimensional faces of $P_X$ and we define
$$W(m,D):= P_X \backslash Z(m,D)$$
the complementary open subset of $P_X$.
Let $$j:W(m,D)\hookrightarrow P_X$$ be the open inclusion of $W(m,D)$ into $P_X$.
\newtheorem{contractible2}[poslem]{Lemma}
\begin{contractible2}\label{cntct2}
For any $\sigma \in \Delta$, $V_\sigma \cap Z(m,D)$ is either empty or it is contractible.
\end{contractible2}
\begin{proof}: Suppose $V_\sigma \cap Z(m,D)= \bigcup_{ \{e \in I \mid \langle m , e \rangle < -d_e \} } (F_e \cap V_\sigma) $ is non-empty. Then each non-empty set $F_e \cap V_\sigma$ in the union contains $T_\sigma$, since 
$$ F_e \cap V_\sigma = \bigcup_{\tau \ni e} T_\tau \cap \bigcup_{\delta \subseteq \sigma } T_\delta $$ and the $\{T_\alpha\}_{\alpha \in \Delta}$ partition $P_X$. They are also convex as $F_e$ and $V_\sigma$ are both convex subsets. If $s \in T_\sigma$, then the constant map $f:V_\sigma \cap Z(m,D) \rightarrow \{ s\}$ is a homotopy equivalence since $f \circ \iota = id_{\{ s\}}$ where $\iota:\{ s\} \hookrightarrow V_\sigma \cap Z(m,D)$ is the inclusion and 
\begin{align*}
 F:(V_\sigma \cap Z(m,D))\! \times \![0,1] &\rightarrow  (V_\sigma \cap Z(m,D)) \\
 (a,t) &\mapsto tm + (1-t)a
\end{align*}
is a well defined, continuous map such that $F(-,0)=id_{V_\sigma \cap Z(m,D)}$ and $F(-,1)= \iota \circ f$ so $\iota \circ f \simeq id_{V_\sigma \cap Z(m,D)}$. Thus $V_\sigma \cap Z(m,D) \simeq \{ s\}$. \end{proof}

\section{Proof of Theorem~\ref{maint} }
We begin with a lemma.
\newtheorem{acyclic}{Lemma}[section]
\begin{acyclic}\label{acycl}
Let $\sigma \in \Delta$ and $j_\sigma : V_\sigma \hookrightarrow P_X$ be the inclusion. Then:
$$ H^i({V_\sigma}, j_\sigma^*j_! \underline{\C}_{W} ) = 0 \quad \forall \; i>0 $$
\end{acyclic}
\begin{proof} First note that ${j_\sigma}^{-1}(W)= W \cap V_\sigma $ and we have the commutative diagram
\[ \begin{CD}
W \cap V_\sigma        @>{j^\sigma}>>           V_\sigma    \\
@V{j_\sigma}VV                               @VV{j_\sigma}V       \\
  W                    @>j>>                     P_X
\end{CD} \]
where $j^\sigma: W \cap V_\sigma \hookrightarrow V_\sigma$ is the inclusion. Then using equation 6.13 on page 111 of Iversen \cite{Iversen}
 $${(j_\sigma)}^* j_!\underline{\C}_{W} = {j^\sigma}_!{(j_\sigma)}^*\underline{\C}_{W} $$ We note that ${(j_\sigma)}^*\underline{\C}_{W}= \underline{\C}_{W \cap V_\sigma} = {(j^\sigma)}^*\underline{\C}_{V_\sigma}$ so 
${(j_\sigma)}^* j_!\underline{\C}_{W} = {j^\sigma}_!{(j^\sigma)}^*\underline{\C}_{V_\sigma}$. Since $j^\sigma$ is the inclusion of an open set in $V_\sigma$ and letting $i^\sigma : V_\sigma \cap Z \hookrightarrow V_\sigma$ denote the inclusion of its complement, there is a well known short exact sequence of sheaves (see~\cite{Hartshorne,Iversen}):
\begin{equation} \label{seseq} 0 \longrightarrow {j^\sigma}_!{(j^\sigma)}^*\underline{\C}_{V_\sigma} \longrightarrow \underline{\C}_{V_\sigma} \longrightarrow {i^\sigma}_*{(i^\sigma)}^*\underline{\C}_{V_\sigma} \longrightarrow 0
\end{equation}
This induces a long exact sequence of cohomology:
\begin{equation*} \label{leseq} \dots \longrightarrow H^i(V_\sigma \cap Z, \C) \longrightarrow H^{i+1}(V_\sigma, {j^\sigma}_!{(j^\sigma)}^*\underline{\C}_{V_\sigma}) \longrightarrow H^{i+1}(V_\sigma, \C) \longrightarrow \dots
\end{equation*}
Using this and applying Lemmas \ref{cntct} and \ref{cntct2} the result follows. \end{proof}

The rest of the proof follows the same strategy as the proof of Demazure's Theorem given in \cite{Dani}.
We start by considering the left hand side of equation (\ref{mthm}). There is a natural covering of $X(\Delta)$ by affine open sets $U_{\sigma}$ with $\sigma \in \Delta$ and intersections of such sets are of the same form. By Serre's Theorem this covering is acyclic, and thus the cohomology $H^i(X, \Oo(D))$ is the same as the $i$-dimensional cohomology of the $\check{\text{C}}$ech complex with this covering:
\begin{equation}\label{cech1} C^*(\{ U_{\sigma}\}_{\sigma \in \Delta}, \Oo(D) ) = ( \dots \stackrel{d}{\longrightarrow} \bigoplus_{\sigma} H^0(U_{\sigma}, \Oo(D) ) \stackrel{d}{\longrightarrow} \cdots)  \end{equation}

There is a natural $M$-grading on each term of the complex and this grading is preserved by the differentials. The $m$-th piece of the cohomology, $ H^i(X,\Oo(D))_m $ equals the $i$-dimensional cohomology of the complex of $m$-th pieces.

Now we look at the right hand side of (\ref{mthm}). From Lemma~\ref{ocover} there is an open covering of $P_X$ by the sets $\{ V_\sigma\}_{\sigma \in \Delta}$, and by Lemma~\ref{acycl} this cover is acyclic. Therefore by Leray's theorem, the cohomology $H^i(P_X,j_! \underline{\C}_{W})$ is the $i$-dimensional cohomology of the $\check{\text{C}}$ech complex
\begin{equation}\label{cech2} C^*(\{ V_\sigma\}_{\sigma \in \Delta}, j_! \underline{\C}_{W} ) = ( \dots \stackrel{d}{\longrightarrow} \bigoplus_{\sigma} H^0({V_\sigma}, j_! \underline{\C}_{W} ) \stackrel{d}{\longrightarrow} \cdots)  \end{equation}
Finally we show that the two spaces $H^0(U_{\sigma}, \Oo(D) )_m$ and $H^0({V_\sigma}, j_! \underline{\C}_{W} )$ are both isomorphic to either $\C$ or $0$ for given $\sigma$ and $m$. Then since the open covers for the two $\check{\text{C}}$ech complexes are indexed by the same set and the differentials are defined in terms of restriction maps which in both cases correspond to the identity map on $\C$ or zero, it can be seen that the cohomology of the complexes is the same.

We know (from \cite{Dani, Fulton}) that $H^0(U_{\sigma}, \Oo(D) )_m = \C$ when $m$ belongs to the set $\{u \in M_{\R} \mid u \geq \psi_D \; \text{on}\; |\sigma| \}$ and is $0$ otherwise. On the other hand consider again the first part of the long exact sequence induced by equation (\ref{seseq}):
$$0 \longrightarrow H^0({V_\sigma},j_{\sigma}^* j_! \underline{\C}_{W} ) \longrightarrow H^0(V_{\sigma}, \C )\longrightarrow H^0({Z(m,D)}\cap{V_\sigma}, \C ) \longrightarrow$$
Obviously $H^0({V_\sigma}, j_! \underline{\C}_{W} )$ is isomorphic to $\C$ when ${Z(m,D)}\cap{V_\sigma}$ is empty and is $0$ otherwise. However 
\begin{align*}
{Z(m,D)}\cap{V_\sigma} = \varnothing 
&\iff F_{e} \cap V_\sigma = \bigcup_{\tau \ni e} T_\tau \cap \bigcup_{\delta \subseteq \sigma } T_\delta = \varnothing \quad \forall {e} \in \{e_i \mid \langle m , e_i \rangle < -d_i \} \\
&\iff \langle m , e_i \rangle \geq -d_i \quad \forall e_i \in |\sigma| \\
&\iff m \in \{u \in M_{\R} \mid u \geq \psi_D \; \text{on}\; |\sigma| \}
\end{align*}
{\it Remark}: We can extend the result to cases when the toric variety is quasi-projective. Suppose $Y=Y(\Sigma)$ is a quasi-projective toric variety embedded via a toric morphism as an open subset of a projective toric variety $X$. There is a natural cover of $Y$ by affine open pieces $\{ U_\sigma \mid \sigma \in \Sigma\}$. On the polytope $P_X$ there is a corresponding open cover $\{ V_\sigma \mid \sigma \in \Sigma\}$ of an open subset $P_Y$ of $P_X$. Restricting to $Y$ and $P_Y$, the rest of the proof follows through; these covers are both acyclic as before, there are two corresponding $\check{\text{C}}$ech complexes (which are sub-complexes of the graded version of (\ref{cech1}), and (\ref{cech2})) and by the computation above it can be seen that they are the same.

\section{Example}\label{example}
As a proof of concept we do a calculation. Any toric quiver variety comes with a natural collection of line bundles, the universal bundles. It was shown by Altmann and Hille (Theorem~3.6 of \cite{AltHille}) that for any smooth Fano toric quiver variety the universal bundles form a strongly exceptional collection. The example here is not Fano and was constructed so that neither the Kawamata-Viehweg vanishing theorem nor Theorem~2.4(ii) from \cite{Mustata} is sufficient 
 to prove the vanishing of the Ext's between all the universal bundles. Thus these vanishing theorems can not be used to show that the collection of universal bundles is strongly exceptional. We show here by direct calculation that it is, and furthermore extend the collection to produce a full strongly exceptional collection.

Consider the following quiver $Q$ with weights at the vertices as labelled.
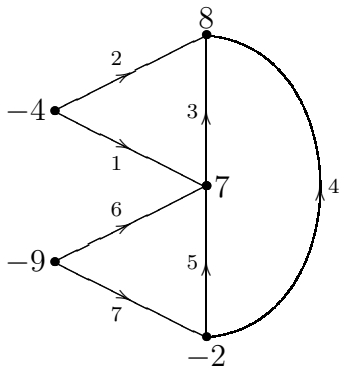
\[ \begin{xy} <10mm,0mm>:
,0*@{*}*^+!U{-2};(0,2)*@{*}**@{-}**@{-}?(0.5)*@{>}*^+!R{\scriptstyle{5}}
,0*@{*};(-2,1)*@{*}**@{-}**@{-}?(0.5)*@{<}*^+!UR{\scriptstyle{7}}
,(-2,1)*@{*};(0,2)*@{*}**@{-}**@{-}?(0.5)*@{>}*^+!DR{\scriptstyle{6}}
,(0,2)*@{*};(0,4)*@{*}**@{-}**@{-}?(0.5)*@{>}*^+!R{\scriptstyle{3}}
,(0,2)*@{*};(-2,3)*@{*}**@{-}**@{-}?(0.5)*@{<}*^+!UR{\scriptstyle{1}}
,(-2,3)*@{*};(0,4)*@{*}**@{-}**@{-}?(0.5)*@{>}*^+!DR{\scriptstyle{2}}
,0;(0,4)**\crv{(2,0.15)&(2,3.85)}?(0.5)*@{>}*^+!L{\scriptstyle{4}}
,(0,2)*@{*}*^+!L{7}
,(-2,1)*@{*}*^+!R{-9}
,(-2,3)*@{*}*^+!R{-4}
,(0,4)*@{*}*^+!D{8}
\end{xy} \]
Let $Q_0$ be the set of vertices and $Q_1$ the set of arrows. There are two maps $h,t:Q_1 \rightarrow Q_0$ taking an arrow to its head and tail respectively. The corresponding toric quiver variety $X$ (see \cite{Hille}), is smooth and complete and has rays generated by $\{ e_1=(e_3 +e_4),e_2=-(e_3 +e_4),e_3,e_4,e_5,e_6=-(e_4+e_5),e_7=(e_4+e_5) \}$. There is a correspondence between the prime torus invariant divisors on $X$ and elements of $Q_1$. The universal bundles $ \{  \Ll_v \mid v \in Q_0 \}$ satisfy the property that $\Ll_{ha} \otimes \Ll_{ta}^* \cong \Oo(E_a)$ for all $a \in Q_1$ and are unique up to a twist. Let $ \{e_3^{\vee},e_4^{\vee},e_5^{\vee}\}$ be the dual basis of $ \{e_3,e_4,e_5\}$ and then $M$ is the integer lattice generated by $ \{e_3^{\vee},e_4^{\vee},e_5^{\vee}\}$. The polytope $P_X$ is three dimensional, the stereographic projection of which is shown below.

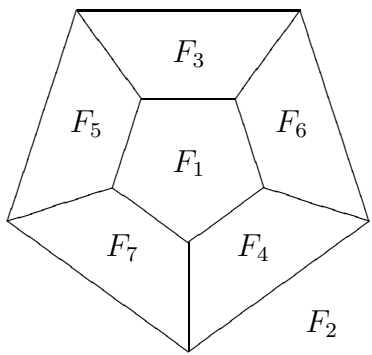
\[\begin{xy}/r25mm/:
,{*{F_1}\xypolygon5"O"{~={198}\dir{}}}
,{\xypolygon5"I"{~:{(0.42,0):}~={198}\dir{}}}
,"O1";"I1"**@{-},"O2";"I2"**@{-},"O3";"I3"**@{-},"O4";"I4"**@{-},"O5";"I5"**@{-}
,(1.7,-0.85)*{F_2},(1,0.58)*{F_3},(0.46,0.21)*{F_5},(1.54,0.21)*{F_6},(0.65,-0.45)*{F_7},(1.34,-0.45)*{F_4}
\end{xy}\]

Consider again the long exact sequence of cohomology induced by the short exact sequence (\ref{seseq3}):
\begin{equation*} \label{long2} \dots \longrightarrow H^i(Z(m,D), \C) \longrightarrow H^{i+1}(P_X,j_! \underline{\C}_{W}) \longrightarrow H^{i+1}(P_X, \C) \longrightarrow \dots
\end{equation*}
where $Z(m,D)$ is some collection of closed faces of the polytope $P_X$. Looking at all possible unions of closed faces we see that topologically, there are five possibilities. From the long exact sequence above it can be seen which of these contributes to each non-zero cohomology, and the table below lists all of these. Thus $H^i(P_X,j_! \underline{\C}_{W}) \cong \C$ for all $Z$ in the row $H^i \cong \C$ and $H^i(P_X,j_! \underline{\C}_{W}) = 0$ for all other $Z$. We use the notation $$F_J:= \bigcup_{j\in J}F_j$$ and let
\begin{align*}
{\bf J}:=\big{\lbrace}&\{1, 2\},\{3, 4\},\{3, 7\},\{4, 5\},\{5, 6\},\{6, 7\},\\
&\{3, 4, 7\},\{3, 4, 5\},\{4, 5, 6\},\{5, 6, 7\},\{3, 6, 7\}\big{\rbrace} \end{align*} 
contain those sets $J$ such that $F_J$ has two connected components. We denote the complement of $J$ in $\{1,2,3,4,5,6,7\}$ by $J^c$. 
 
\[ \begin{tabular}{|m{3cm}|m{4cm}|m{4.05cm}|}
\hline $H^0 \cong \C$ & $Z(m,D) \simeq \varnothing$ & $Z=\varnothing$ \\
\hline $H^1 \cong \C$&$Z(m,D) \simeq \{pt\} \times \{ pt\}$& $ Z = F_J$ for  $J \in \bf{J} $ \\ 
\hline $H^2 \cong \C$&$Z(m,D) \simeq S^1$& $Z=F_J $ for $J^c \in \bf{J}$  \\ 
\hline $H^3 \cong \C$&$Z(m,D) \simeq S^2$& $Z=P_X$  \\ 
\hline
\end{tabular} \]

Let $\Oo(D)$ be any invertible line bundle on $X$. We can choose $D$ to be of the form $D:= d_1E_1 + d_2 E_2 + d_6 E_6 + d_7 E_7$, a torus-invariant Weil divisor 
which is Cartier, since $\Pic(X)$ is rank 4 and generated by $\{E_1, E_2, E_6, E_7\}$. Then from the table above, using the decomposition in equation~(\ref{dsd}) and applying Theorem~\ref{maint}, it can be seen that $ H^0(X,\Oo(D)) \neq 0$ if and only if there exists $m \in M$ such that $Z(m,D) = \varnothing$. This holds if and only if
\begin{align*} & \exists m:=m_3{e_3}^{\vee}+m_4{e_4}^{\vee}+m_5{e_5}^{\vee} \in M \; \text{such that}\; \langle m , e_i \rangle \geq -d_i \; \forall i \\
\iff &\exists m_1,m_2,m_3 \in \Z \; \text{where}\; m_i \geq 0, -d_1 \leq m_3+m_4 \leq d_2, -d_7 \leq m_4+m_5 \leq d_6 \\
\iff &0 \leq (d_1 + d_2),  0 \leq (d_6 + d_7),d_2 \geq 0, d_6 \geq 0
\end{align*}
Similarly $ H^1(X,\Oo(D)) \neq 0$ if and only if there exists $m \in M$ such that $Z(m,D)=F_J $ for some $J^c \in \bf{J}$. Treating each case separately, we obtain a list of all the regions in $\Pic(X)$ where $ H^1(X,\Oo(D))$ is nonzero. We can then produce a complete list of regions in $\Pic(X)$ where each cohomology group is nonzero either by continuing the process or applying Serre Duality. The list of these 24 regions is given in tabular form below.
\[ \begin{tabular}{|m{1.5cm}|m{0.4cm}|m{12cm}|}
\hline $ H^0 \neq 0$& 1 & $0 \leq (d_1 + d_2),\;  0 \leq (d_6 + d_7), \; 0 \leq d_2,\; 0 \leq d_6$\\
\hline $H^1 \neq 0$& 2 &$-2 \geq (d_1 + d_2),\;  0 \leq (d_6 + d_7),\;-1 \geq d_1,\; 0 \leq d_6$  \\ 
& 3 & $0 \leq (d_1 + d_2),\;  0 \leq (d_6 + d_7),\;2 \leq d_1,\; 1 \leq (d_1+d_6)$  \\
& 4 & $0 \leq (d_1 + d_2),\;  0 \leq d_6,\;-1 \geq d_7,\; 1 \leq (d_1+d_6),\; 2 \leq (d_1-d_7)$  \\
& 5 & $0 \leq (d_1 + d_2),\;  0 \leq (d_6 + d_7),\;2 \leq d_7,\; 1 \leq (d_2+d_7)$  \\
& 6 & $0 \leq (d_1 + d_2),\;  1 \leq (d_2 + d_7),\;0 \leq d_2,\; 2 \leq (d_2-d_6)$  \\
& 7 & $0 \leq (d_1 + d_2),\;  -2 \geq (d_6 + d_7),\;0 \leq d_2,\; -1 \geq d_7$  \\
& 8 & $0 \leq (d_1 + d_2),\;  1 \leq (d_1 + d_6),\;2 \leq d_1,\; 2 \leq (d_1-d_7)$  \\
& 9 & $0 \leq (d_1 + d_2),\;  0 \leq (d_6 + d_7),\;2 \leq d_1,\; 2 \leq d_7$  \\
& 10 & $0 \leq (d_1 + d_2),\;  1 \leq (d_2 + d_7),\;-3 \geq d_6,\;2 \leq d_7,\; 2 \leq (d_2-d_6)$  \\
& 11 & $0 \leq (d_1 + d_2),\;  -2 \geq (d_6 + d_7),\;0 \leq d_2,\; 2 \leq (d_2-d_6)$  \\
& 12 & $0 \leq (d_1 + d_2),\;  -2 \geq (d_6 + d_7),\;-1 \geq d_7,\; 2 \leq (d_1-d_7)$ \\ 
\hline $H^2 \neq 0$& 13 & $0 \leq (d_1 + d_2),\;  -2 \geq (d_6 + d_7),\;2 \leq d_1,\; -3 \geq d_6$\\ 
& 14 & $-2 \geq (d_1 + d_2),\;  -2 \geq (d_6 + d_7),\;-1 \geq d_1,\; -3 \geq (d_1+d_6)$  \\
& 15 & $-2 \geq (d_1 + d_2),\;  -3 \geq d_6,\;2 \leq d_7,\; -3 \geq (d_1+d_6),\; -2 \geq (d_1-d_7)$  \\
& 16 & $-2 \geq (d_1 + d_2),\;  -2 \geq (d_6 + d_7),\;-1 \geq d_7,\; -3 \geq (d_2+d_7)$  \\
& 17 & $-2 \geq (d_1 + d_2),\;  -3 \geq (d_2 + d_7),\;-3 \geq d_2,\; -2 \geq (d_2-d_6)$  \\
& 18 & $-2 \geq (d_1 + d_2),\;  0 \leq (d_6 + d_7),\;-3 \geq d_2,\; 2 \leq d_7$  \\
& 19 & $-2 \geq (d_1 + d_2),\;  -3 \geq (d_1 + d_6),\;-1 \geq d_1,\; -2 \geq (d_1-d_7)$  \\
& 20 & $-2 \geq (d_1 + d_2),\;  -2 \geq (d_6 + d_7),\;-1 \geq d_1,\; -1 \geq d_7$  \\
& 21 & $-2 \geq (d_1 + d_2),\;  -3 \geq (d_2 + d_7),\;0 \leq d_6,\;-1 \geq d_7,\; -2 \geq (d_2-d_6)$  \\
& 22 & $-2 \geq (d_1 + d_2),\;  0 \leq (d_6 + d_7),\;-3 \geq d_2,\; -2 \geq (d_2-d_6)$  \\
& 23 & $-2 \geq (d_1 + d_2),\;  0 \leq (d_6 + d_7),\;2 \leq d_7,\; -2 \geq (d_1-d_7)$ \\ \hline 
$H^3 \neq 0$& 24 & $-2 \geq (d_1 + d_2),\;  -2 \geq (d_6 + d_7),\; -3 \geq d_2,\; -3 \geq d_6$ \\ \hline
\end{tabular} \]
Using this table it can easily be confirmed that the universal bundles form an exceptional collection. For each pair of line bundles $\Ll_p $ and $\Ll_q$ in the collection we wish to calculate the cohomology of $\Ll_q \otimes \Ll_p^* \cong \Oo(D)$, where $D$ is chosen  so that it is a linear combination of the generators $\{E_1, E_2, E_6, E_7\}$ of Pic(X). For example, if $p$ is the vertex with weight 8 and $q$ is the vertex with weight -2, then $\Ll_q \otimes \Ll_p^* \cong \Oo(-E_1 + E_2 + E_6 - E_7)$. Looking at the table, the only region which contains the point $(d_1,d_2,d_6,d_7)=(-1,1,1,-1)$ is region 1 and so there is no higher cohomology. Similarly there is no higher cohomology for all other pairs of line bundles in the universal collection. Thus the Ext's between all pairs of line bundles in the collection are zero and the universal bundles form a strongly exceptional collection.

\subsection{Buchsbaum-Rim}
We digress slightly to describe a generalised Koszul complex called the Buchsbaum-Rim complex, see \cite{Northcott} Appendix~C. We will then use this in our example to construct a collection of line bundles which extends the collection of universal bundles $ \{  \Ll_v \mid v \in Q_0 \}  $ on $X$ and which by construction spans the derived category $D^b(X)$. Then using the cohomology calculation above once more, we show that this collection is actually strongly exceptional.

Let $V$,$W$ be vector bundles of ranks $m$,$n$ over an arbitrary base, and let $f:W \rightarrow V$ be a bundle map. Then the Buchsbaum-Rim complex $(K_*,d_*)$ is as follows:
$K_0=V$, $K_1=W$ with $d_1=f$ and then
$$ K_{r+1}=\Lambda^{m+r+1}(W)\otimes S^r(V^*) \otimes \operatorname{det}V^*  $$
For $n>2$, the maps $d_n$ are defined to be interior product with $f:W \rightarrow V$ regarded as a section of $W^* \otimes V$ and $d_2:\Lambda^{m+1}(W) \otimes \operatorname{det}V^* \rightarrow W$ is interior product with $ \Lambda^m(f): \Lambda^m(W) \rightarrow \Lambda^m(V) $ regarded as a section of $\Lambda^m(W^*) \otimes \Lambda^m(V)$.
The complex $(K_*,d_*)$ is exact away from the support of coker$f$, and in particular since it has length $n-m+1$, if the support of coker$f$ is in codimension $n-m+1$, then the complex is a resolution of coker$f$ (\cite{BuchsEisen}).

Any toric quiver variety $X$ comes with a presentation of the diagonal in $X \times X$:
$$ \bigoplus_{a \in Q_1}\Ll_{ta}\boxtimes \Ll^*_{ha} \longrightarrow  \bigoplus_{v \in Q_0}\Ll_{v}\boxtimes \Ll^*_{v} $$
where the components of this map are $\phi_a \boxtimes 1 - 1 \boxtimes \phi_a^*$. The ranks of these bundles are $n=|Q_1|$ and $m=|Q_0|$ respectively, and the support of the cokernel of this map (i.e. the diagonal in $X \times X$) has codimension equal to the dimension of $X$, namely $n-m+1$. Therefore from above we can see that Buchsbaum-Rim complex for this presentation is a resolution of the diagonal $\Oo_\Delta$ in $X \times X$. Each term of the complex consists of a product of line bundles on each side of the $\boxtimes$. Taking all the line bundles that appear on either one side or the other gives a collection which spans $D^b(X)$ (\cite{King}). We have:
\begin{align*} S^r(V^*) \otimes \operatorname{det}(V^*) &= \bigoplus_{\genfrac{}{}{0pt}{}{P:Q_0 \rightarrow \{1,2,\dots\}}{|P|=m+r}}  \left( \bigotimes_{v \in Q_0} \Ll_v^{-P(v)} \right) \boxtimes \left( \bigotimes_{v \in Q_0} \Ll_v^{P(v)} \right)\\
\Lambda^{m+r+1}(W) &= \bigoplus_{\genfrac{}{}{0pt}{}{R:Q_1 \rightarrow \{0,1\}}{|R|=m+r+1}}  \left( \bigotimes_{a \in Q_1} \Ll_{ta}^{R(a)} \right) \boxtimes \left( \bigotimes_{a \in Q_1} \Ll_{ha}^{-R(a)} \right)
\end{align*}
where $ |P|:= \sum_{v \in Q_0} P(v)$ and $ |R|:= \sum_{a \in Q_1} R(a)$. 
Therefore the duals of the line bundles that appear on the right hand side of the $\boxtimes$ are the universal bundles $\Ll_v$, for $v \in Q_0$ and 
$$ \bigotimes_{a \in Q_1} \Ll_{ha}^{R(a)} \otimes \bigotimes_{v \in Q_0} \Ll_v^{-P(v)} $$
for every $P:Q_0 \rightarrow \{ 1,2, \dots\}$ and $R:Q_1 \rightarrow \{0,1\}$ such that $m+1 \leq |P|+1=|R| \leq n$.

Calculating these in the example above, we obtain the collection of ten line bundles corresponding to the vertices in the quiver below. The arrows which correspond to Hom's between the bundles are decorated with the labels of the corresponding divisors, i.e. an arrow from $\Ll_p $ to $\Ll_q$ is labelled 1+5+7 when $\Ll_q \otimes \Ll_p^* \cong \Oo(E_1+E_5+E_7)$.

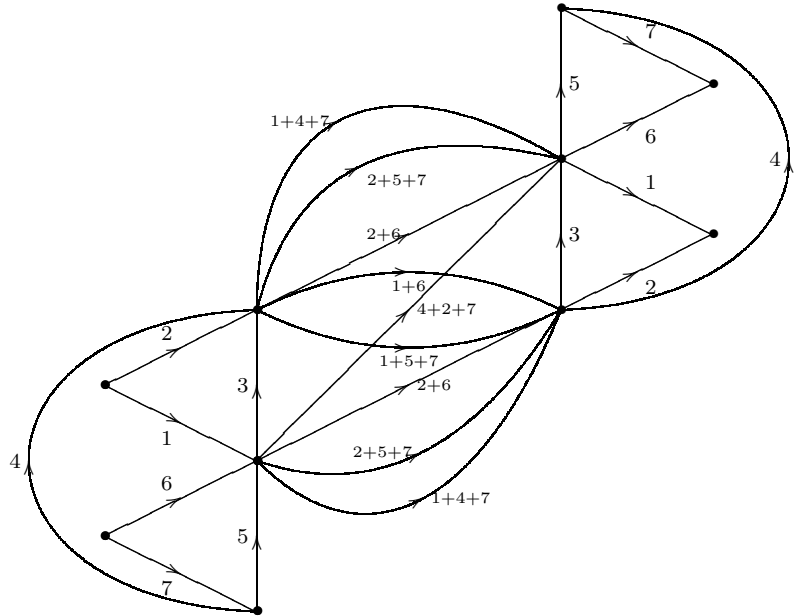
\[\begin{xy} <10mm,0mm>:
,0*@{*};(0,2)*@{*}**@{-}**@{-}?(0.5)*@{>}*^+!R{\scriptstyle{5}}
,0*@{*};(-2,1)*@{*}**@{-}**@{-}?(0.5)*@{<}*^+!UR{\scriptstyle{7}}
,(-2,1)*@{*};(0,2)*@{*}**@{-}**@{-}?(0.5)*@{>}*^+!DR{\scriptstyle{6}}
,(0,2)*@{*};(0,4)*@{*}**@{-}**@{-}?(0.5)*@{>}*^+!R{\scriptstyle{3}}
,(0,2)*@{*};(-2,3)*@{*}**@{-}**@{-}?(0.5)*@{<}*^+!UR{\scriptstyle{1}}
,(-2,3)*@{*};(0,4)*@{*}**@{-}**@{-}?(0.5)*@{>}*^+!DR{\scriptstyle{2}}
,(0,4);(4,4)**\crv{(2,5)}?(0.5)*@{>}*^+!U{\scriptscriptstyle{1+6}}
,(0,4);(4,4)**\crv{(2,3)}?(0.5)*@{>}*^+!U{\scriptscriptstyle{1+5+7}}
,(4,4)*@{*};(4,6)**@{-}**@{-}?(0.5)*@{>}*^+!L{\scriptstyle{3}}
,(4,4);(6,5)**@{-}**@{-}?(0.5)*@{>}*^+!UL{\scriptstyle{2}}
,(4,6)*@{*};(6,5)*@{*}**@{-}**@{-}?(0.5)*@{>}*^+!DL{\scriptstyle{1}}
,(4,6);(6,7)*@{*}**@{-}**@{-}?(0.5)*@{>}*^+!UL{\scriptstyle{6}}
,(4,6);(4,8)*@{*}**@{-}**@{-}?(0.5)*@{>}*^+!L{\scriptstyle{5}}
,(4,8);(6,7)**@{-}**@{-}?(0.5)*@{>}*^+!DL{\scriptstyle{7}}
,(0,2);(4,4)**@{-}**@{-}?(0.5)*@{>}*^+!L{\scriptscriptstyle{2+6}}
,(0,4);(4,6)**@{-}**@{-}?(0.5)*@{>}*^+!R{\scriptscriptstyle{2+6}}
,(0,2);(4,6)**@{-}**@{-}?(0.5)*@{>}*^+!L{\scriptscriptstyle{4+2+7}}
,0;(0,4)**\crv{(-4,0.15)&(-4,3.85)}?(0.5)*@{>}*^+!R{\scriptstyle{4}}
,(4,4);(4,8)**\crv{(8,4.15)&(8,7.85)}?(0.5)*@{>}*^+!R{\scriptstyle{4}}
,(0,2);(4,4)**\crv{(2.3,1.2)}?(0.5)*@{>}*^+!R{\scriptscriptstyle{2+5+7}}
,(0,2);(4,4)**\crv{(1.5,0.5)&(3,1.5)}?(0.5)*@{>}*^+!L{\scriptscriptstyle{1+4+7}}
,(0,4);(4,6)**\crv{(0.7,6.8)}?(0.5)*@{>}*^+!UL{\scriptscriptstyle{2+5+7}}
,(0,4);(4,6)**\crv{(0,6.5)&(1.5,7.5)}?(0.5)*@{>}*^+!R{\scriptscriptstyle{1+4+7}}
\end{xy}\]
Again for each pair of line bundles $\Ll_p $ and $\Ll_q$ in the collection we wish to calculate the cohomology of $\Ll_q \otimes \Ll_p^* \cong \Oo(D)$. By looking at the table above it can be seen that for each such pair, the divisor $D$ does not lie in any of the regions of Pic$(X)$ where $H^i(\Oo(D)) \neq 0$ for $i >0$. Hence there is no higher cohomology, so the Ext's between all pairs of line bundles in the collection are zero. Thus we have a full strongly exceptional collection of line bundles on our toric variety.

\subsection{Bondal's Collection}
Given any smooth toric variety $X$, Bondal has described a method to produce a candidate collection of line bundles on $X$, which for a certain class of Fano varieties is expected to be strongly exceptional. In this section we determine this collection in the case of our non-Fano example from section~\ref{example}, and show that it is not strongly exceptional. This is not unexpected but gives another illustration of the method.

For any toric variety $X$ and $l \in \N$, there is a well-defined toric morphism 
$$ \pi_l : {X} \rightarrow X$$
which restricts, on the torus $T$, to the map
$$ \pi_l : T \rightarrow T , \quad t \mapsto t^l $$ 
In the case when $X$ is smooth then the direct image,
$$ (\pi_l)_* \mathcal{O}_X = \bigoplus_\chi \mathcal{L}_\chi  $$ is a direct sum of line bundles indexed by the characters of the $l$-torsion subgroup of $T$. This is because the map $\pi_l$ is the quotient of $X$ by this group. 
The set $\B$ of line bundles which occur as summands of this direct sum for all sufficiently large $l$ exists and is given by
\begin{align*}
\B =& \{ \Oo({D}) 
\mid {D}:= -\sum_{i=1}^d \{  \langle e_i , m \rangle  \} E_i,\; m \in M_{\Q} \}\\
=&\{ \Oo(\bar{D}) 
\mid \bar{D}:= \sum_{i=1}^d \left\lfloor  \langle e_i , m \rangle  \right\rfloor E_i,\; m = \sum_{i=1}^n m_ie_1^\vee \in M_{\Q}, \; 0\leq m_i <1 \}
\end{align*}
where $M_\Q = M \otimes_\Z \Q$ and $ \{ \alpha \} = \alpha - \left\lfloor \alpha \right\rfloor \geq 0$ is the fractional part of $\alpha$, for $\alpha \in \mathbb Q$.
(For a more general construction one can instead consider $ (\pi_l)_* \mathcal{O}_X({D}) $ for some divisor ${D}$.)
In the example this produces a collection of 12 line bundles corresponding to the vertices of the quiver below, where the top vertex corresponds to $\Oo_X$:

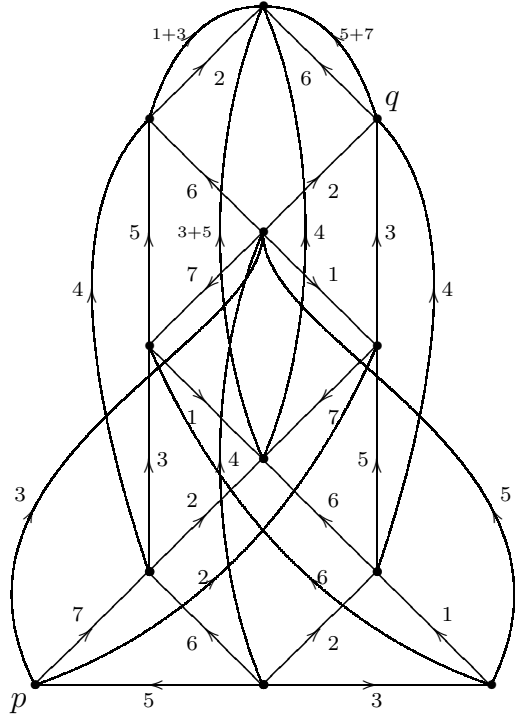
\[\begin{xy} <15mm,0mm>:
,0*@{*}*^+!RU{p};(2,0)*@{*}**@{-}**@{-}?(0.5)*@{<}*^+!U{\scriptstyle{5}}
,0*@{*};(1,1)*@{*}**@{-}**@{-}?(0.5)*@{>}*^+!DR{\scriptstyle{7}}
,(2,0)*@{*};(4,0)*@{*}**@{-}**@{-}?(0.5)*@{>}*^+!U{\scriptstyle{3}}
,(2,0)*@{*};(1,1)*@{*}**@{-}**@{-}?(0.5)*@{>}*^+!UR{\scriptstyle{6}}
,(2,0)*@{*};(3,1)*@{*}**@{-}**@{-}?(0.5)*@{>}*^+!UL{\scriptstyle{2}}
,(1,1)*@{*};(2,2)*@{*}**@{-}**@{-}?(0.5)*@{>}*^+!DR{\scriptstyle{2}}
,(1,1);(1,5)**\crv{(0,4)}?(0.5)*@{>}*^+!R{\scriptstyle{4}}
,(3,1);(3,5)**\crv{(4,4)}?(0.5)*@{>}*^+!L{\scriptstyle{4}}
,(1,1)*@{*};(1,3)**@{-}**@{-}?(0.5)*@{>}*^+!L{\scriptstyle{3}}
,(3,1);(2,2)**@{-}**@{-}?(0.5)*@{>}*^+!DL{\scriptstyle{6}}
,(3,1)*@{*};(3,3)*@{*}**@{-}**@{-}?(0.5)*@{>}*^+!R{\scriptstyle{5}}
,(3,3);(3,5)*@{*}*^+!DL{q}**@{-}**@{-}?(0.5)*@{>}*^+!L{\scriptstyle{3}}
,(3,3);(2,2)*@{*}**@{-}**@{-}?(0.5)*@{>}*^+!UL{\scriptstyle{7}}
,(1,3)*@{*};(2,2)**@{-}**@{-}?(0.5)*@{>}*^+!UR{\scriptstyle{1}}
,(1,3);(1,5)**@{-}**@{-}?(0.5)*@{>}*^+!R{\scriptstyle{5}}
,(2,4);(1,3)**@{-}**@{-}?(0.5)*@{>}*^+!DR{\scriptstyle{7}}
,(2,4);(3,3)**@{-}**@{-}?(0.5)*@{>}*^+!DL{\scriptstyle{1}}
,(4,0);(3,1)**@{-}**@{-}?(0.5)*@{>}*^+!DL{\scriptstyle{1}}
,(2,4)*@{*};(1,5)*@{*}**@{-}**@{-}?(0.5)*@{>}*^+!UR{\scriptstyle{6}}
,(2,4);(3,5)**@{-}**@{-}?(0.5)*@{>}*^+!UL{\scriptstyle{2}}
,(3,5);(2,6)*@{*}**@{-}**@{-}?(0.5)*@{>}*^+!UR{\scriptstyle{6}}
,(1,5);(2,6)**@{-}**@{-}?(0.5)*@{>}*^+!UL{\scriptstyle{2}}
,(2,0);(2,4)**\crv{(1.5,1.15)&(1.5,2.85)}?(0.5)*@{>}*^+!L{\scriptstyle{4}}
,(2,2);(2,6)**\crv{(2.5,3.15)&(2.5,4.85)}?(0.5)*@{>}*^+!L{\scriptstyle{4}}
,(2,2);(2,6)**\crv{(1.5,3.15)&(1.5,4.85)}?(0.5)*@{>}*^+!R{\scriptscriptstyle{3+5}}
,(0,0);(2,4)**\crv{(-1,2)&(2,3)}?(0.3)*@{>}*^+!DR{\scriptstyle{3}}
,(4,0);(2,4)**\crv{(5,2)&(2,3)}?(0.3)*@{>}*^+!DL{\scriptstyle{5}}
,(0,0);(3,3)**\crv{(2,0.7)}?(0.45)*@{>}*^+!R{\scriptstyle{2}}
,(4,0);(1,3)**\crv{(2,0.7)}?(0.45)*@{>}*^+!L{\scriptstyle{6}}
,(1,5);(2,6)**\crv{(1.3,6)}?(0.5)*@{>}*^+!R{\scriptscriptstyle{1+3}}
,(3,5);(2,6)**\crv{(2.7,6)}?(0.5)*@{>}*^+!L{\scriptscriptstyle{5+7}}
\end{xy}\]
Note that the quiver: 
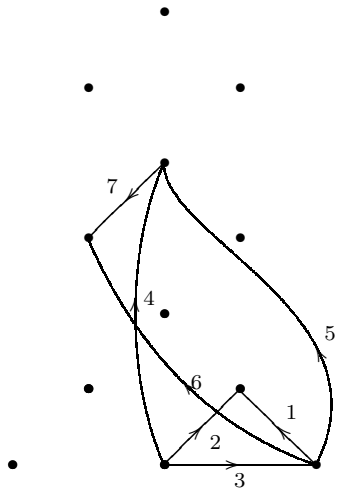
\[\begin{xy} <10mm,0mm>:
,0*@{*};(2,0)*@{*}
,0*@{*};(1,1)*@{*}
,(2,0)*@{*};(4,0)*@{*}**@{-}**@{-}?(0.5)*@{>}*^+!U{\scriptstyle{3}}
,(2,0)*@{*};(1,1)*@{*}
,(2,0)*@{*};(3,1)*@{*}**@{-}**@{-}?(0.5)*@{>}*^+!UL{\scriptstyle{2}}
,(1,1)*@{*};(2,2)*@{*}
,(1,1);(1,5)
,(3,1);(3,5)
,(1,1)*@{*};(1,3)
,(3,1);(2,2)
,(3,1)*@{*};(3,3)*@{*}
,(3,3);(3,5)*@{*}
,(3,3);(2,2)*@{*}
,(1,3)*@{*};(2,2)
,(1,3);(1,5)
,(2,4);(1,3)**@{-}**@{-}?(0.5)*@{>}*^+!DR{\scriptstyle{7}}
,(2,4);(3,3)
,(4,0);(3,1)**@{-}**@{-}?(0.5)*@{>}*^+!DL{\scriptstyle{1}}
,(2,4)*@{*};(1,5)*@{*}
,(2,4);(3,5)
,(3,5);(2,6)*@{*}
,(1,5);(2,6)
,(2,0);(2,4)**\crv{(1.5,1.15)&(1.5,2.85)}?(0.55)*@{>}*^+!L{\scriptstyle{4}}
,(2,2);(2,6)
,(0,0);(2,4)
,(4,0);(2,4)**\crv{(5,2)&(2,3)}?(0.3)*@{>}*^+!DL{\scriptstyle{5}}
,(0,0);(3,3)
,(4,0);(1,3)**\crv{(2,0.7)}?(0.5)*@{>}*^+!L{\scriptstyle{6}}
\end{xy}\]
is a subquiver. This subquiver is also contained in the quiver corresponding to the full strongly exceptional collection on $X$ which we produced above. It can be seen that, up to a twist, the duals of the universal bundles are contained in $ \B $. There are however $\operatorname{Ext}^1$'s between some of the line bundles in the collection. Consider $\Ll_p$ and $\Ll_q$ where $p$ and $q$ are labelled on the quiver above. It can be seen from the quiver that $\Ll_q \otimes \Ll_p^* \cong \Oo(-E_1+2E_2)$, so $d_1=-1, d_2=2 , d_6=d_7=0$. Then looking back at the table from the cohomology calculation, we see that these satisfy the sixth set of inequalities, so 
$$\operatorname{Ext}^1(\Ll_p,\Ll_q) \cong H^1(\Ll_q \otimes \Ll_p^*) \neq 0.$$
We find that there are no $\operatorname{Ext}^2$'s between any of the line bundles. All the $\operatorname{Ext}^1$'s are shown below:

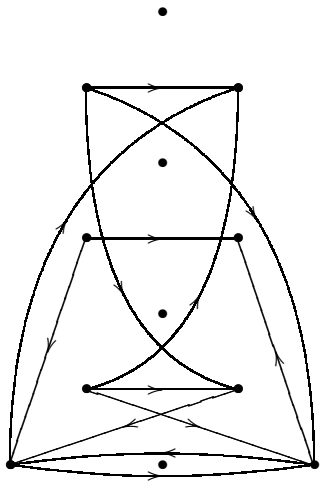
\[\begin{xy} <10mm,0mm>:
,0*@{*};(3,1)*@{*}**@{-}**@{-}?(0.5)*@{<}
,0*@{*};(3,5)*@{*}**\crv{(0,4)}?(0.5)*@{>}
,(1,3)*@{*};(0,0)*@{*}**@{-}**@{-}?(0.5)*@{>}
,(1,1)*@{*};(3,5)*@{*}**\crv{(3,1.5)}?(0.5)*@{>}
,(1,3)*@{*};(3,3)*@{*}**@{-}**@{-}?(0.5)*@{>}
,(1,5)*@{*};(3,5)*@{*}**@{-}**@{-}?(0.5)*@{>}
,(1,1);(3,1)**@{-}**@{-}?(0.5)*@{>}
,(1,1)*@{*};(4,0)*@{*}**@{-}**@{-}?(0.5)*@{>}
,(4,0);(3,3)*@{*}**@{-}**@{-}?(0.5)*@{>}
,(4,0);(1,5)*@{*}**\crv{(4,4)}?(0.5)*@{<}
,(3,1)*@{*};(1,5)**\crv{(1,1.5)}?(0.5)*@{<}
,(1,3);(2,0)*@{*}
,(2,2)*@{*};(1,3)
,(2,4);(3,3)
,(4,0);(3,1)
,(2,4)*@{*};(1,5)*@{*}
,(2,4);(3,5)
,(3,5);(2,6)*@{*}
,(1,5);(2,6)
,(0,0);(4,0)**\crv{(1.5,-0.2)&(2.5,-0.2)}?(0.5)*@{>}
,(4,0);(0,0)**\crv{(2.5,0.2)&(1.5,0.2)}?(0.5)*@{>}
\end{xy}\]
One may check that no subcollection containing any ten of these line bundles is strongly exceptional and thus there is no full strongly exceptional subcollection of $\B$.

Mathematical Sciences,\\
University of Bath,\\
Bath,\\
BA2 7AY\\
U.K.\\
\url{ntb20@maths.bath.ac.uk}
\end{document}